\documentclass[graybox]{svmult}

\usepackage{mathptmx}       % selects Times Roman as basic font
\usepackage{helvet}         % selects Helvetica as sans-serif font
\usepackage{courier}        % selects Courier as typewriter font
\usepackage{type1cm}        % activate if the above 3 fonts are
\usepackage{makeidx}         % allows index generation
\usepackage{graphicx}        % standard LaTeX graphics tool
\usepackage{multicol}        % used for the two-column index
\usepackage[bottom]{footmisc}% places footnotes at page bottom

\usepackage{amsmath,amssymb,amsopn,bm}
\usepackage{url}
\usepackage{anslistings} % source code listings
\usepackage{shortvrb}
\usepackage{stmaryrd} % extra symbols
\usepackage{graphicx}
\usepackage[caption=false]{subfig}
\usepackage{booktabs}
\usepackage{tcolorbox}
\tcbuselibrary{skins,breakable}

\newtcolorbox{mytcolorbox}{%
colback=white,
boxrule=0.2mm,
arc=1pt,
left=0mm,right=0mm,top=0mm,bottom=0mm,middle=0mm
}

\usepackage[colorinlistoftodos,prependcaption,textsize=footnotesize]{todonotes}
\reversemarginpar

\usepackage[normalem]{ulem}

\makeatletter
\g@addto@macro\bfseries{\boldmath}
\makeatother

\newcommand{\R}{\mathbb{R}}

\DeclareMathOperator*{\divv}{div}
\DeclareMathOperator*{\diam}{diam}

\newcommand{\bfe}{\bm{e}}
\newcommand{\bfn}{\bm{n}}

\newcommand{\bfx}{\bm{x}}

\newcommand{\bfp}{\bm{p}}

\newcommand{\bfv}{\bm{v}}
\newcommand{\bfu}{\bm{u}}

\newcommand{\bff}{\bm{f}}

\newcommand{\bfzero}{\bm{0}}

\newcommand{\bfV}{\bm{V\!}} % NB: make a little tighter
\newcommand{\bfW}{\bm{W\!}} % NB: make a little tighter

\newcommand{\mcK}{\mathcal{K}}

\newcommand{\mcO}{\mathcal{O}}

\newcommand{\hatOmega}{\widehat{\Omega}}
\newcommand{\hatmcK}{\widehat{\mcK}}

\newcommand{\OO}{\mathcal{O}}

\makeindex             % used for the subject index
\begin{document}

\title*{A MultiMesh Finite Element Method for the Stokes Problem}

  %% High order methods on arbitrarily many intersecting meshes:
  %% Multimesh}
%\titlerunning{Multimesh}
%\author{Name of First Author and Name of Second Author}
\author{August Johansson \and
  %Benjamin Kehlet \and
  Mats G. Larson \and
  Anders Logg}
%\institute{Name of First Author \at Name, Address of Institute, \email{name@email.address}
%\and Name of Second Author \at Name, Address of Institute \email{name@email.address}}
\institute{August Johansson \at Simula Research Laboratory, P.O.\ Box 134, 1325 Lysaker, Norway, \email{august.johansson@gmail.com}
 % \and
 % Benjamin Kehlet \at  Simula Research Laboratory, P.O.\ Box 134, 1325 Lysaker, Norway, \email{benjamik@simula.no}
  \and
  Mats G. Larson \at Department of Mathematics, Ume{\aa} University, 90187, Ume{\aa}, Sweden, \email{mats.larson@math.umu.se}
  \and
  Anders Logg \at Department of Mathematical Sciences, Chalmers University of Technology and University of Gothenburg, 41296 G\"oteborg, Sweden, \email{logg@chalmers.se}.}

\maketitle

%%%%%%%%%%%%%%%%%%%%%%%%%%%%%%%%%%%%%%%%%%%%%%%%%%%%%%%%%%%%%%%%%%%%%%%%%%%%%%%%
\abstract{The multimesh finite element method enables the solution of partial differential equations on a computational mesh composed by multiple arbitrarily overlapping meshes. The discretization is based on a continuous--discontinuous function space with interface conditions enforced by means of Nitsche's method.
In this contribution, we consider the Stokes problem as a first step towards flow applications. The multimesh formulation leads to so called cut elements in the underlying meshes close to overlaps. These demand stabilization to ensure coercivity and stability of the stiffness matrix. We employ a consistent least-squares term on the overlap to ensure that the inf-sup condition holds. We here present the method for the Stokes problem, discuss the implementation, and verify that we have optimal convergence.
%Each chapter should be preceded by an abstract (10--15 lines long) that summarizes the content. The abstract will appear \textit{online} at \url{www.SpringerLink.com} and be available with unrestricted access. This allows unregistered users to read the abstract as a teaser for the complete chapter. As a general rule the abstracts will not appear in the printed version of your book unless it is the style of your particular book or that of the series to which your book belongs.\newline\indent
%  Please use the 'starred' version of the new Springer \texttt{abstract} command for typesetting the text of the online abstracts (cf. source file of this chapter template \texttt{abstract}) and include them with the source files of your manuscript. Use the plain \texttt{abstract} command if the abstract is also to appear in the printed version of the book.
  }
%%%%%%%%%%%%%%%%%%%%%%%%%%%%%%%%%%%%%%%%%%%%%%%%%%%%%%%%%%%%%%%%%%%%%%%%%%%%%%%%
\section{Introduction}
\label{sec:1}

Consider the Stokes problem
\begin{alignat}{2}
  \label{eq:stokes1}
  -\Delta \bfu + \nabla p &= \bff \qquad &&\text{in $\Omega$},
  \\
  \label{eq:stokes2}
  \divv \bfu &= 0 \qquad &&\text{in $\Omega$},
  \\
  \bfu &= \bfzero \qquad &&\text{on $\partial \Omega$},
\end{alignat}
for the velocity $\bfu : \Omega \rightarrow \R^d$ and pressure $\bfp : \Omega \rightarrow \R$ in a polytopic domain $\Omega\subset\R^d$, $d = 2,3$.

The Stokes problem is considered here as a first step towards a multimesh formulation for multi-body flow problems, and ultimately fluid--structure interaction problems, in which each body is discretized by an individual boundary-fitted mesh and the boundary-fitted meshes move freely on top of a fixed background mesh. The applications for such a formulation are many, e.g., the simulation of blood platelets in a blood stream, the optimization of the configuration of an array of wind turbines, or the investigation of the effect of building locations in a simulation of urban wind conditions and pollution. Common to these applications is that the multimesh method removes the need for costly mesh (re)generation and allows the platelets, wind turbines or buildings to be moved around freely in the domain, either in each timestep as a part of a dynamic simulation, or in each iteration as part of an optimization problem.

The multimesh formulation presented here is a generalization of the formulation presented and analyzed in~\cite{Johansson:2015aa} for two domains. For comparison, the multimesh discretization of the Poisson problem for arbitrarily many intersecting meshes is presented in~\cite{mmfem-1} and analyzed in~\cite{mmfem-2}.

%%%%%%%%%%%%%%%%%%%%%%%%%%%%%%%%%%%%%%%%%%%%%%%%%%%%%%%%%%%%%%%%%%%%%%%%%%%%%%%%
\section{Notation}
\label{sec:notation}

We first review the notation for domains, interfaces, meshes and overlaps used to formulate the multimesh finite element method. For a more detailed exposition, we refer to~\cite{mmfem-1}.

\begin{figure}
  \centering
  \subfloat[]{\label{fig:three_domains_a}\includegraphics[height=0.2\linewidth]{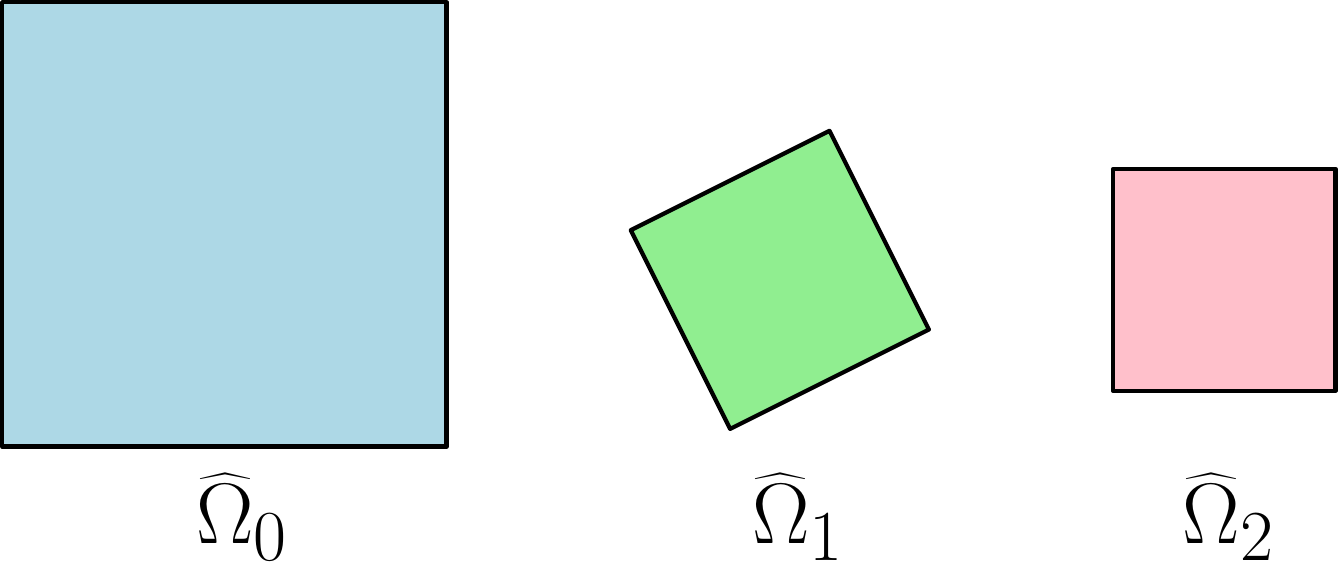}}\qquad\qquad\qquad
  \subfloat[]{\label{fig:three_domains_b}\includegraphics[height=0.2\linewidth]{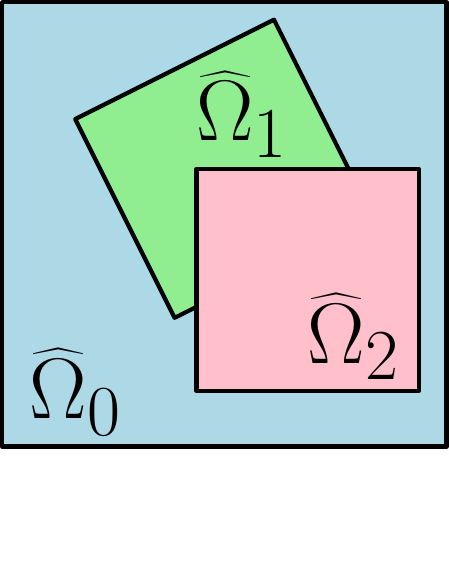}}
  \caption{(a) Three polygonal predomains. (b) The predomains are placed on top of each other in an ordering such that
    $\hatOmega_0$ is placed lowest, $\hatOmega_1$ is in the middle and $\hatOmega_2$ is on top.}
  \label{fig:three_domains}
\end{figure}

\begin{figure}
  \begin{center}
    \includegraphics[height=0.2\linewidth]{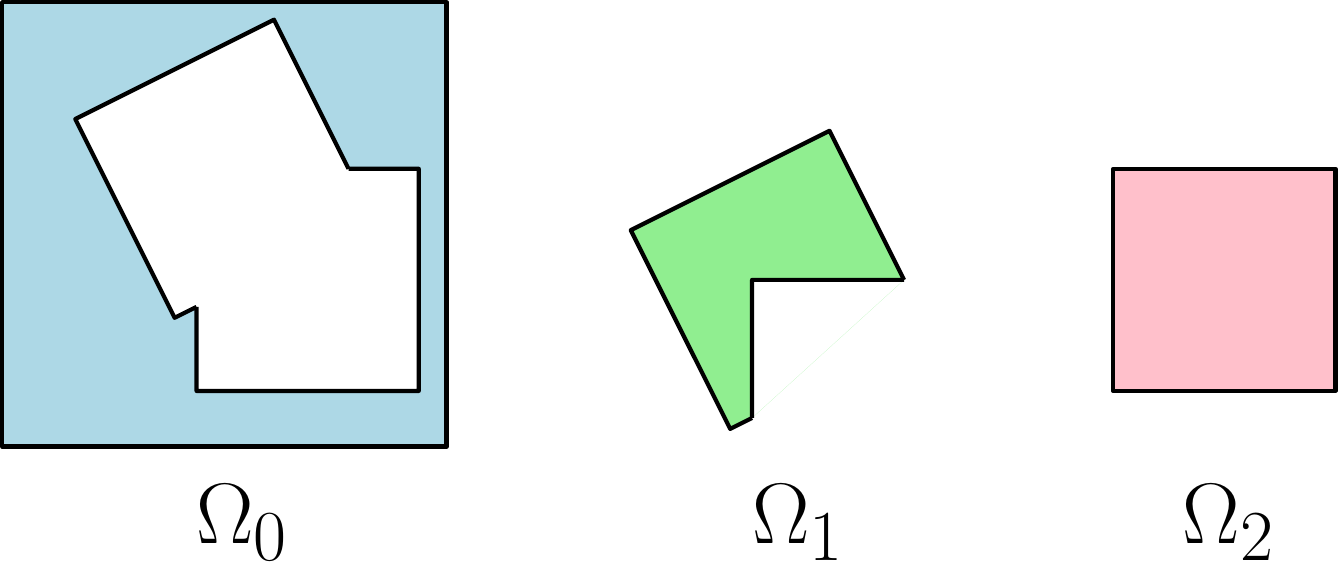}
    \caption{Partition of $\Omega = \Omega_0 \cup \Omega_1 \cup \Omega_2$. Note that $\Omega_2 = \hatOmega_2$.}
    \label{fig:three_domains_partition}
  \end{center}
\end{figure}

\begin{mytcolorbox}
\emph{Notation for domains}
\tcblower

Let $\Omega = \hatOmega_0 \subset \R^d$, $d = 2,3$, be a domain with polytopic boundary (the background domain).

Let $\hatOmega_i \subset \hatOmega_0$, $i=1,\ldots, N$ be the so-called \emph{predomains} with polytopic boundaries (see Figure~\ref{fig:three_domains}).

Let $\Omega_i = \hatOmega_i \setminus \bigcup_{j=i+1}^{N} \hatOmega_j$,  $i=0, \ldots, N$ be a partition of $\Omega$ (see Figure~\ref{fig:three_domains_partition}).
\end{mytcolorbox}

\begin{remark}
  \label{rem:boundary-overlap}
To simplify the presentation, the domains $\Omega_1, \ldots, \Omega_N$ are not allowed to intersect the boundary of $\Omega$.
\end{remark}

\begin{figure}
  \centering
  \subfloat[]{\label{fig:two_interfaces_a}\includegraphics[width=0.2\linewidth]{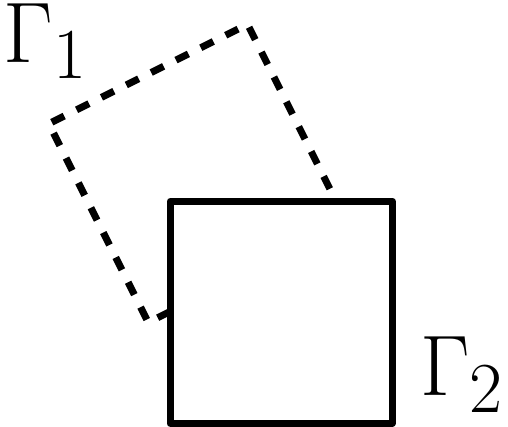}}\qquad\qquad\qquad
  \subfloat[]{\label{fig:two_interfaces_b}\includegraphics[width=0.2\linewidth]{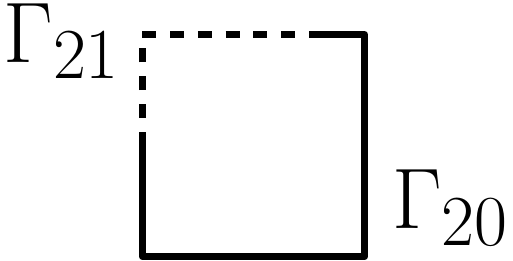}}
  \caption{(a) The two interfaces of the domains in Figure~\ref{fig:three_domains}: $\Gamma_1 = \partial \hatOmega_1 \setminus \hatOmega_2$ (dashed line) and  $\Gamma_2 = \partial \hatOmega_2$ (filled line). Note that $\Gamma_1$ is not a closed curve. (b) Partition of $\Gamma_2 = \Gamma_{20} \cup \Gamma_{21}$.}
\end{figure}

\begin{mytcolorbox}
\emph{Notation for interfaces}
\tcblower

Let the \emph{interface} $\Gamma_i$ be defined by $\Gamma_i = \partial \hatOmega_i \setminus \bigcup_{j=i+1}^N \hatOmega_j$, $i=1, \ldots, N-1$ (see Figure~\ref{fig:two_interfaces_a}).

Let $\Gamma_{ij} = \Gamma_i \cap \Omega_j$, $i > j$ be a partition of $\Gamma_i$ (see Figure~\ref{fig:two_interfaces_b}).
\end{mytcolorbox}

\begin{figure}
  \centering
  \subfloat[]{\label{fig:three_meshes_a}\includegraphics[height=0.2\linewidth]{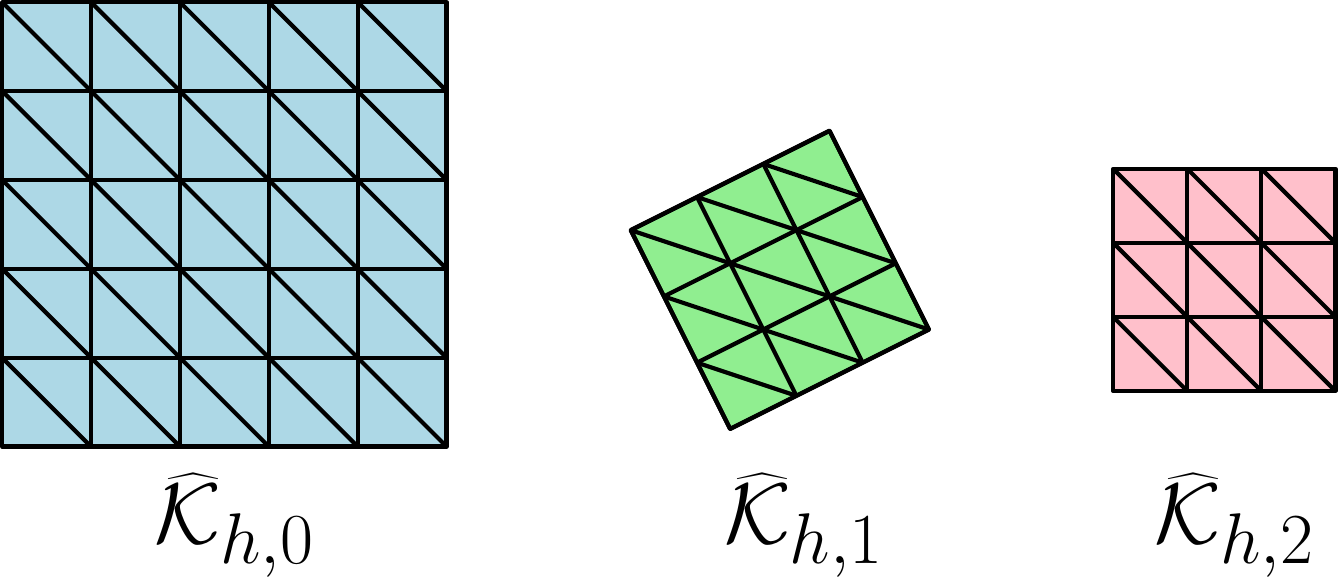}}\qquad\qquad\qquad
  \subfloat[]{\label{fig:three_meshes_b}\includegraphics[height=0.2\linewidth]{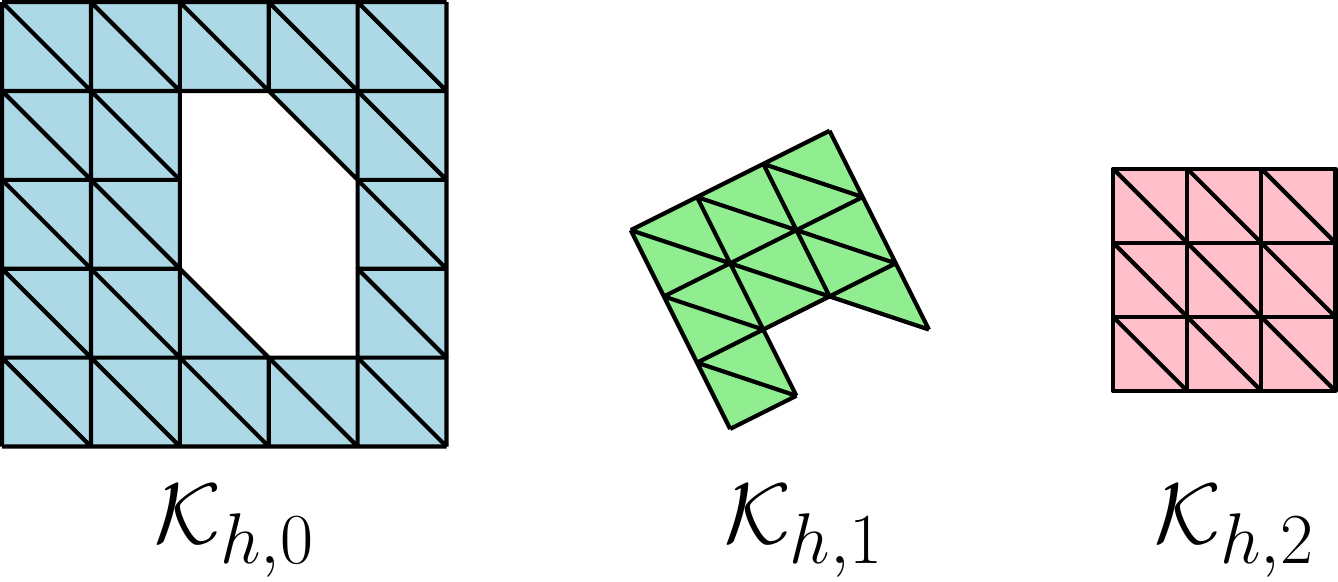}}
  \caption{(a) The three premeshes. (b) The corresponding active meshes (cf.\ Figure~\ref{fig:three_domains}).}
\end{figure}

\begin{figure}
  \centering
  \subfloat[]{\label{fig:overlap}\includegraphics[height=0.2\linewidth]{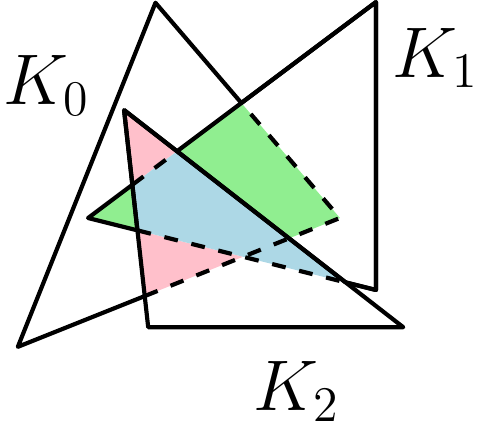}}\qquad\qquad\qquad
  \subfloat[]{\label{fig:multimesh}\includegraphics[height=0.2\linewidth]{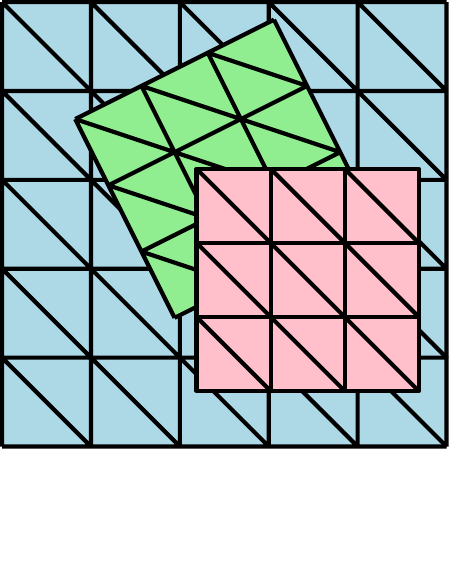}}
  \caption{(a) Given three ordered triangles $K_0$, $K_1$ and $K_2$, the overlaps are $\mcO_{01}$ in green, $\mcO_{02}$ in red and $\mcO_{12}$ in blue. (b) The multimesh of the domains in Figure~\ref{fig:three_domains_b} consists of the active meshes in Figure~\ref{fig:three_meshes_b}.}
\end{figure}

\begin{mytcolorbox}
\emph{Notation for meshes}
\tcblower

Let $\hatmcK_{h,i}$ be a quasi-uniform \cite{BreSco08} \emph{premesh} on $\hatOmega_i$ with mesh parameter $h_i = \max_{K\in \hatmcK_{h,i}} \diam(K)$, $i=0,\ldots,N$ (see Figure~\ref{fig:three_meshes_a}).

Let $\mcK_{h,i} = \{ K \in \hatmcK_{h,i} : K \cap \Omega_i \neq \emptyset \}$, $i=0,\ldots,N$ be the \emph{active meshes} (see Figure~\ref{fig:three_meshes_b}).

The \emph{multimesh} is formed by the active meshes placed in the given ordering (see Figure~\ref{fig:multimesh}).

Let $\Omega_{h,i} = \bigcup_{K\in\mcK_{h,i}} K$, $i=0,\ldots,N$ be the \emph{active domains}.
\end{mytcolorbox}

\begin{mytcolorbox}
\emph{Notation for overlaps}
\tcblower

Let $\OO_i$ denote the \emph{overlap} defined by  $\OO_i = \Omega_{h,i} \setminus \Omega_i$, $i=0,\ldots,N-1$.

Let $\OO_{ij} = \OO_i \cap \Omega_j = \Omega_{h,i} \cap \Omega_j$, $i < j$ be a partition of $\OO_i$.

% For $i < j$, let
%   \begin{align}
%     \label{eq:indicatorfunction}
%     \delta_{ij} =
%     \begin{cases}
%       1, \quad \OO_{ij} \neq \emptyset,
%       \\
%       0, \quad \text{otherwise},
%     \end{cases}
%   \end{align}
% be a function indicating which overlaps are non-empty.
% For ease of notation, we further let $\delta_{ii} = 1$ for $i = 0, \ldots,N$.

% Let $N_{\OO} = \max(\max_i \sum_j \delta_{ij}, \max_j \sum_i \delta_{ij})$ be the maximum number of overlaps. Note that $N_\OO$ is bounded by $N$ but is usually much smaller.
\end{mytcolorbox}

%%%%%%%%%%%%%%%%%%%%%%%%%%%%%%%%%%%%%%%%%%%%%%%%%%%%%%%%%%%%%%%%%%%%%%%%%%%%%%%%
\section{MultiMesh Finite Element Method}

To formulate the multimesh finite element for the Stokes problem~\eqref{eq:stokes1}--\eqref{eq:stokes2}, we assume for each (active) mesh $\mcK_{h,i}$ the existence of a pair of inf-sub stable spaces $\bfV_{h,i}\times Q_{h,i}$, $i=0,1, \ldots, N$, away from the interface. To be precise, we assume inf-sup stability in $\omega_{h,i} \subset \Omega_{h,i}$, where $\omega_{h,i}$ is close to $\Omega_i$ in the sense that $\Omega_{h,i} \setminus \omega_{h,i} \subset U_\delta(\Gamma_i)$, where
\begin{align}
  U_\delta(\Gamma_i) = \bigcup_{\bfx \in \Gamma_i} B_\delta(\bfx)
\end{align}
and $B_\delta(\bfx)$ is a ball of radius $\delta$ centered at $\bfx$. In other words, $U_\delta(\Gamma_i)$ is the tubular neighborhood of $\Gamma$ with thickness $\delta$. In the numerical examples, we let $\omega_{h,i}$ be the union of elements in $\mcK_{h,i}$ with empty intersection with $\Gamma_{ij}$, $j > i$.

The inf-sup condition may expressed on each submesh $\omega_{h,i}$ by
\begin{align}
  \label{assum:infsup}
  \| p_i - \lambda_{\omega_{h,i}}(p) \|_{\omega_{h,i}}
  \lesssim
  \sup_{\bfv \in \bfW_{h,i}} \frac{(\divv \bfv,p)_{\omega_{h,i}}}{\|D \bfv \|_{\omega_{h,i}}},
\end{align}
where $\lambda_{\omega_{h,i}}(p)$ is the average of $p$ over $\omega_{h,i}$ and $\bfW_{h,i}$ is the subspace of $\bfV_{h,i}$ defined by
\begin{align}
  \label{eq:bfW}
  \bfW_{h,i}
  &=
  \{\bfv \in \bfV_{h,i}: \text{$\bfv=\bfzero$ \text{on} $\overline{\Omega_{h,i} \setminus \omega_{h,i}}$}\}.
\end{align}

We now define the multimesh finite element space as the direct sum
\begin{align}
  \bfV_h \times Q_h = \bigoplus_{i=0}^N \bfV_{h,i}\times Q_{h,i},
\end{align}
where $\bfV_h$ and $Q_h$ consist of piecewise polynomial of degree $k$ and $l$, respectively. This means that an element $\bfv \in \bfV_h$ is a tuple $(\bfv_0, \ldots, \bfv_N)$, and the inclusion $\bfV_h \hookrightarrow L^2(\Omega)$ is defined by $\bfv(\bfx) = \bfv_i(\bfx)$ for $\bfx\in\Omega_i$. A similar interpretation is done for $q \in Q_h$.

We now consider the following asymmetric finite element method: Find $(\bfu_h,p_h) \in \bfV_h \times Q_h$ such that
$A_h((\bfu_h, p_h), (\bfv, q)) = l_h(\bfv)$ for all
$(\bfv,q) \in \bfV_h \times Q_h$,
where
\begin{align}
  \label{eq:Ah}
  A_h((\bfu, p), (\bfv, q))
  &=
  a_h(\bfu,\bfv) + s_h(\bfu, \bfv) +  b_h(\bfu,q)
  + b_h(\bfv,p) + d_h((\bfu,p),(\bfv,q)),
  \\
  \label{eq:ah}
  a_{h}(\bfu,\bfv) &= \sum_{i=0}^N (D \bfu_i, D \bfv_i)_{\Omega_i}
  \\ \nonumber
  &\qquad
  - \sum_{i=1}^N \sum_{j=0}^{i-1} \big( (\langle (D \bfu) \cdot \bfn_i \rangle,[ \bfv])_{\Gamma_{ij}}
  + ([ \bfu], \langle (D \bfv) \cdot \bfn_i \rangle)_{\Gamma_{ij}} \big)
  \\ \nonumber
  &\qquad
  +  \sum_{i=1}^N \sum_{j=0}^{i-1} \beta_0 h^{-1}([\bfu],[ \bfv])_{\Gamma_{ij}},
  \\
  s_{h}(\bfu,\bfv)&=
  \sum_{i=0}^{N-1} \sum_{j=i+1}^N  \beta_1 ([D\bfu_i], [D\bfv_i])_{\OO_{ij}},
  \\
  \label{eq:bh}
  b_h(\bfu,q) &= -  \sum_{i=0}^N (\divv \bfu_i,q_i)_{\Omega_i}
  +  \sum_{i=1}^N \sum_{j=0}^{i-1} ([\bfn_i \cdot \bfu],\langle q \rangle )_{\Gamma_{ij}},
  \\
  \label{eq:dh}
  d_h((\bfu,p),(\bfv,q))&=
  \sum_{i=0}^N \delta h^2(\Delta \bfu_i - \nabla p_i,\Delta \bfv_i + \nabla q_i)_{\Omega_{h,i}\setminus \omega_{h,i}},
  \\
  l_h(\bfv) &=  \sum_{i=0}^N (\bff,\bfv_i)_{\Omega_i}
  -  \sum_{i=0}^N \delta h^2(\bff,\Delta \bfv_i + \nabla q_i)_{\Omega_{h,i}\setminus \omega_{h,i}}.
\end{align}
Here, $\beta_0$ and $\beta_1$ are stabilization parameters that must be sufficiently large to ensure that the bilinear form $A_h$ is coercive; cf.~\cite{Johansson:2015aa} for an analysis of the two-domain case.

For simplicity, we use the global mesh size $h$ here and throughout the presentation. If the meshes are of substantially different sizes, it may be beneficial to introduce the individual mesh sizes $h_i$ in~\eqref{eq:dh} and the average $h_i^{-1} + h_j^{-1}$ in~\eqref{eq:ah}.

Note that since $\Gamma_i$ is partitioned into interfaces $\Gamma_{ij}$ relative to underlying meshes, the sums of the interface terms are over  $0 \leq j < i \leq N$. In contrast, the sums of the overlap terms are over $0 \leq i < j \leq N$ since the overlap $\OO_i$ is partitioned into overlaps $\OO_{ij}$ relative to overlapping meshes.

The jump terms on $\OO_{ij}$ and $\Gamma_{ij}$ are defined by $[\bfv] = \bfv_i - \bfv_j$, where $\bfv_i$ and $\bfv_j$ are the finite element solutions (components) on the active meshes $\mcK_{h,i}$ and $\mcK_{h,j}$. The average normal flux is defined on $\Gamma_{ij}$ by
\begin{equation}\label{eq:average}
  \langle \bfn_i\cdot \nabla \bfv \rangle = (\bfn_i \cdot \nabla \bfv_{i} + \bfn_{i} \cdot \nabla \bfv_{j})/2.
\end{equation}
Here, any convex combination is valid~\cite{Hansbo:2003aa}.

The proposed formulation~\eqref{eq:Ah} is identical to the one proposed in \cite{Johansson:2015aa} with sums over all domains and interfaces. Also note the similarity with the multimesh formulation for the Poisson problem presented in~\cite{mmfem-1}, the difference being the additional least-squares term $d_h$ (and the corresponding term in $l_h$) since we only assume inf-sup stability in $\omega_{h,i}$. If we do not assume inf-sup stability anywhere (e.g.\ if we would use a velocity-pressure element of equal order), the least-squares term should be applied over the whole domain as in \cite{Massing:2014aa}. Please cf.\ \cite{Massing:2014aa} for the use of a symmetric $d_h$.

Other stabilization terms may be considered. By norm equivalence, the stabilization term $s_h(\bfu, \bfv)$ may alternatively be formulated as
\begin{equation}\label{eq:sh-L2}
  s_h (\bfu,\bfv) = \sum_{i=0}^{N-1} \sum_{j=i+1}^N \beta_2 h^{-2} ([ \bfu ], [\bfv ])_{\OO_{ij}}.
\end{equation}
where $\beta_2$ is a stabilization parameter; see~\cite{mmfem-2}.

Note that the finite element method weakly approximates continuity in the sense that $[\bfu_h] = \bfzero$ and $[\bfn_i \cdot \nabla \bfu_h] = 0$ on all interfaces.

%% Stabilization of the pressure may also be desired since the method only guarantees $p_h \in L^2(\Omega)$, and thus one may consider using
%% \begin{align}
%%   \sum_{i=0}^{N-1} \sum_{j=i+1}^N \beta_3 h^{-2} ([p],[q])_{\OO_{ij}},
%% \end{align}
%% where $\beta_3$ is a stabilization parameter.

%%%%%%%%%%%%%%%%%%%%%%%%%%%%%%%%%%%%%%%%%%%%%%%%%%%%%%%%%%%%%%%%%%%%%%%%%%%%%%%%
\section{Implementation}

We have implemented the multimesh finite element method as part of the software framework FEniCS~\cite{Logg:2012aa,Alnaes:2015aa}. One of the main features of FEniCS is the form language UFL~\cite{Alnaes:2014aa} which allows variational forms to be expressed in near-mathematical notation. However, to express the multimesh finite element method~\eqref{eq:Ah}, a number of custom measures must be introduced. In particular, new measures must be introduced for integrals over cut cells, interfaces and overlaps.
These measures are then mapped to quadrature rules that are computed at runtime. An overview of these algorithms algorithms and the implementation is given in~\cite{Johansson:2017ab}.

To express the multimesh finite element method, we let \texttt{dX} denote the integration over domains $\Omega_i$, $i=0, \ldots, N$, including cut cells. Integration over $\Gamma_{ij}$ and $\OO_{ij}$ are expressed using the measures \texttt{dI} and \texttt{dO}, respectively. We let \texttt{dC} denote integration over $\Omega_{h,i} \setminus \omega_{h,i}$. Now the multimesh finite element method for the Stokes problem may be expressed as
\begin{lstlisting}[language=Python,numbers=none]
  a_h = inner(grad(u), grad(v))*dX \
      - inner(avg(grad(u)), tensor_jump(v, n))*dI \
      - inner(avg(grad(v)), tensor_jump(u, n))*dI \
      + beta_0/h * inner(jump(u), jump(v))*dI
  s_h = beta_1 * inner(jump(grad(u)), jump(grad(v)))*dO
  b_h = lambda v, q: inner(-div(v), q)*dX \
                   + inner(jump(v, n), avg(q))*dI
  d_h = delta*h**2 * inner(-div(grad(u)) + grad(p), \
                           -div(grad(v)) - grad(q))*dC
\end{lstlisting}
This makes it easy to implement the somewhat lengthy form~\eqref{eq:Ah}, as well as investigate the effect of different stabilization terms.

%%%%%%%%%%%%%%%%%%%%%%%%%%%%%%%%%%%%%%%%%%%%%%%%%%%%%%%%%%%%%%%%%%%%%%%%%%%%%%%%
\section{Numerical Results}

\begin{figure}
  \begin{center}
    \includegraphics[width=0.320\linewidth]{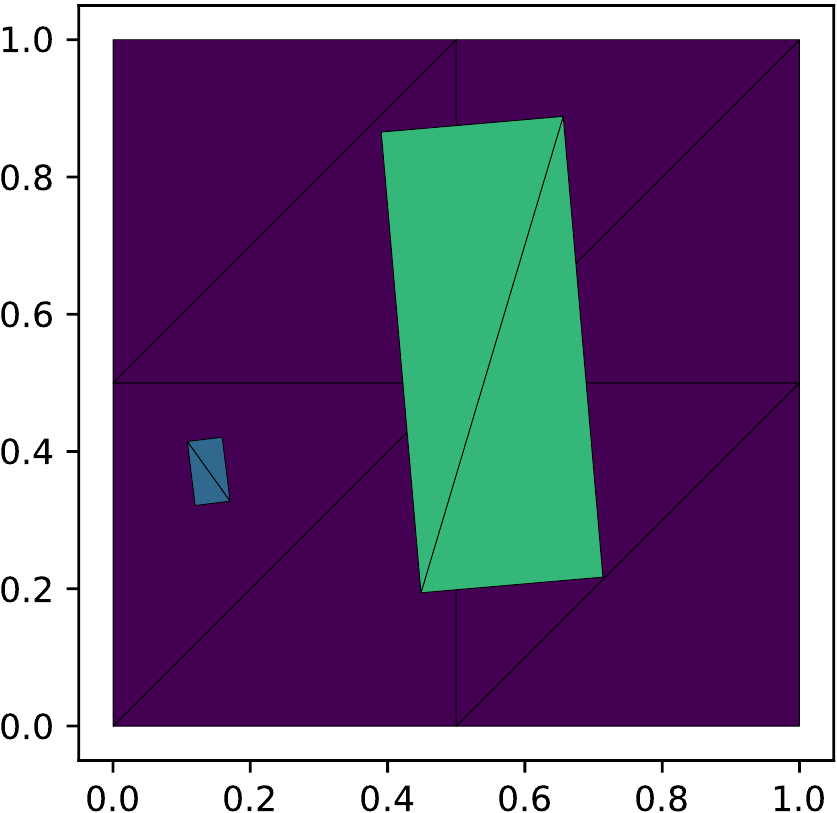}\ \
    \includegraphics[width=0.320\linewidth]{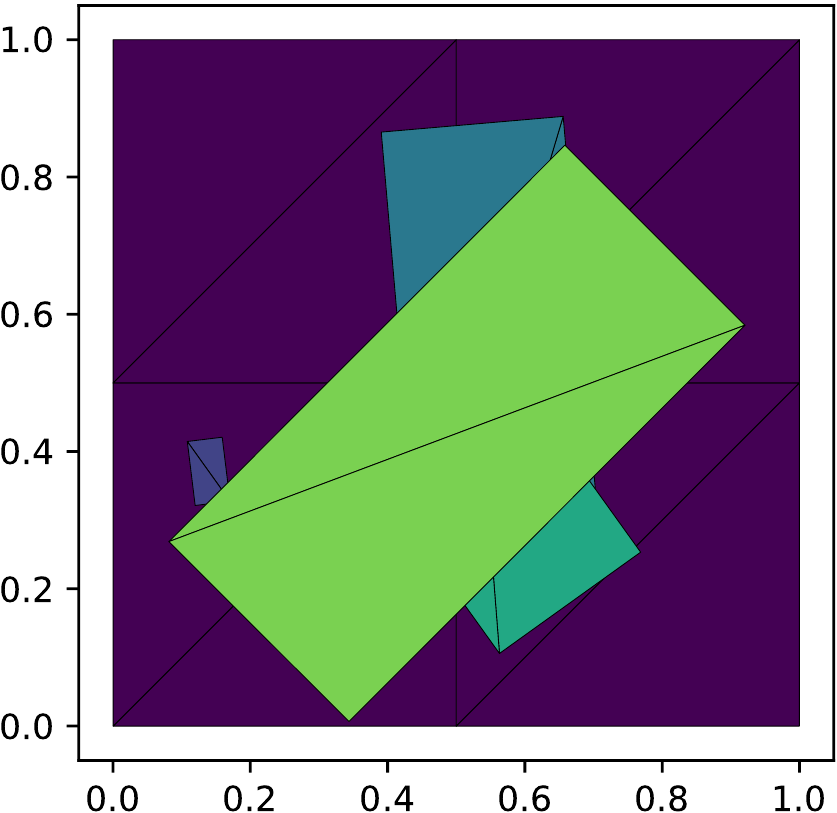}\ \ %\\ \ \\
    \includegraphics[width=0.320\linewidth]{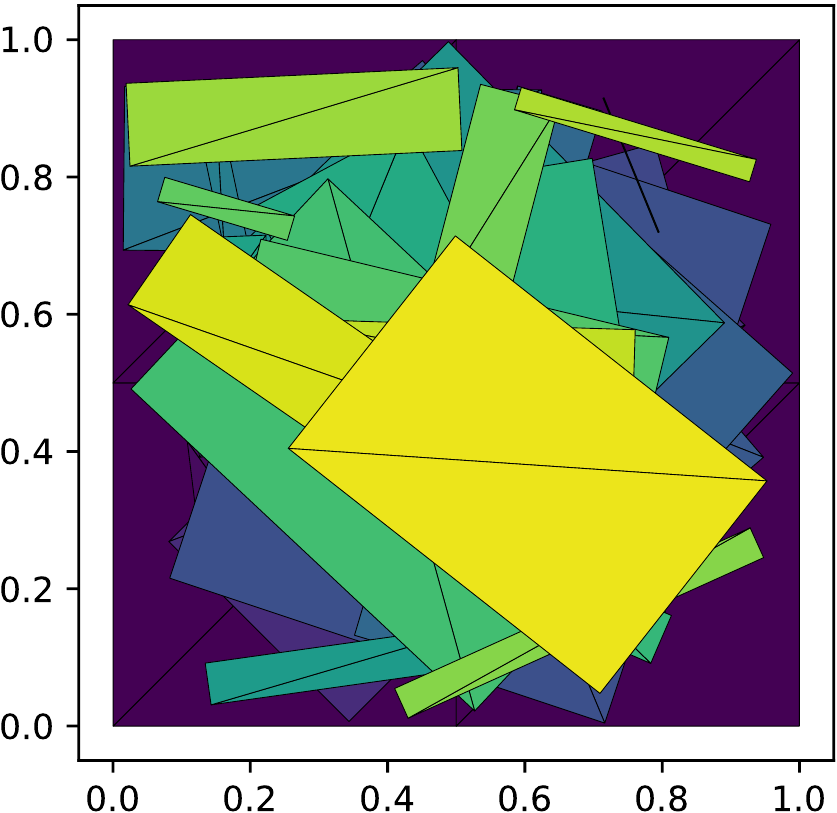}
  \end{center}
  \caption{A sequence of $N$ meshes are randomly placed on top of a fixed background mesh of the unit square shown here for $N=2, 4$ and $32$ using the coarsest refinement level.}
  \label{fig:poisson_meshes}
\end{figure}

\begin{figure}
  \centering

  \includegraphics[width=0.32\textwidth]{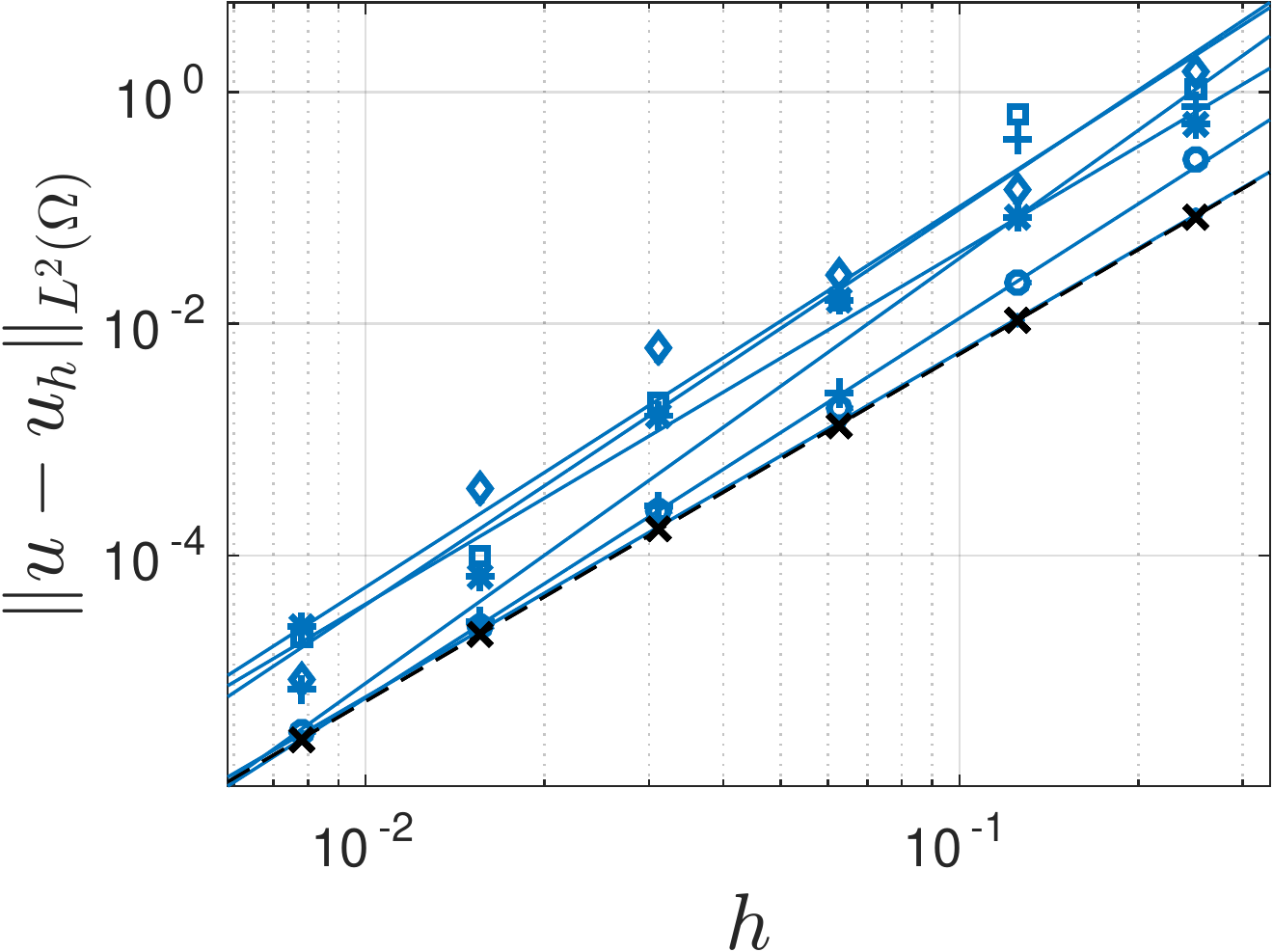}
  \includegraphics[width=0.32\textwidth]{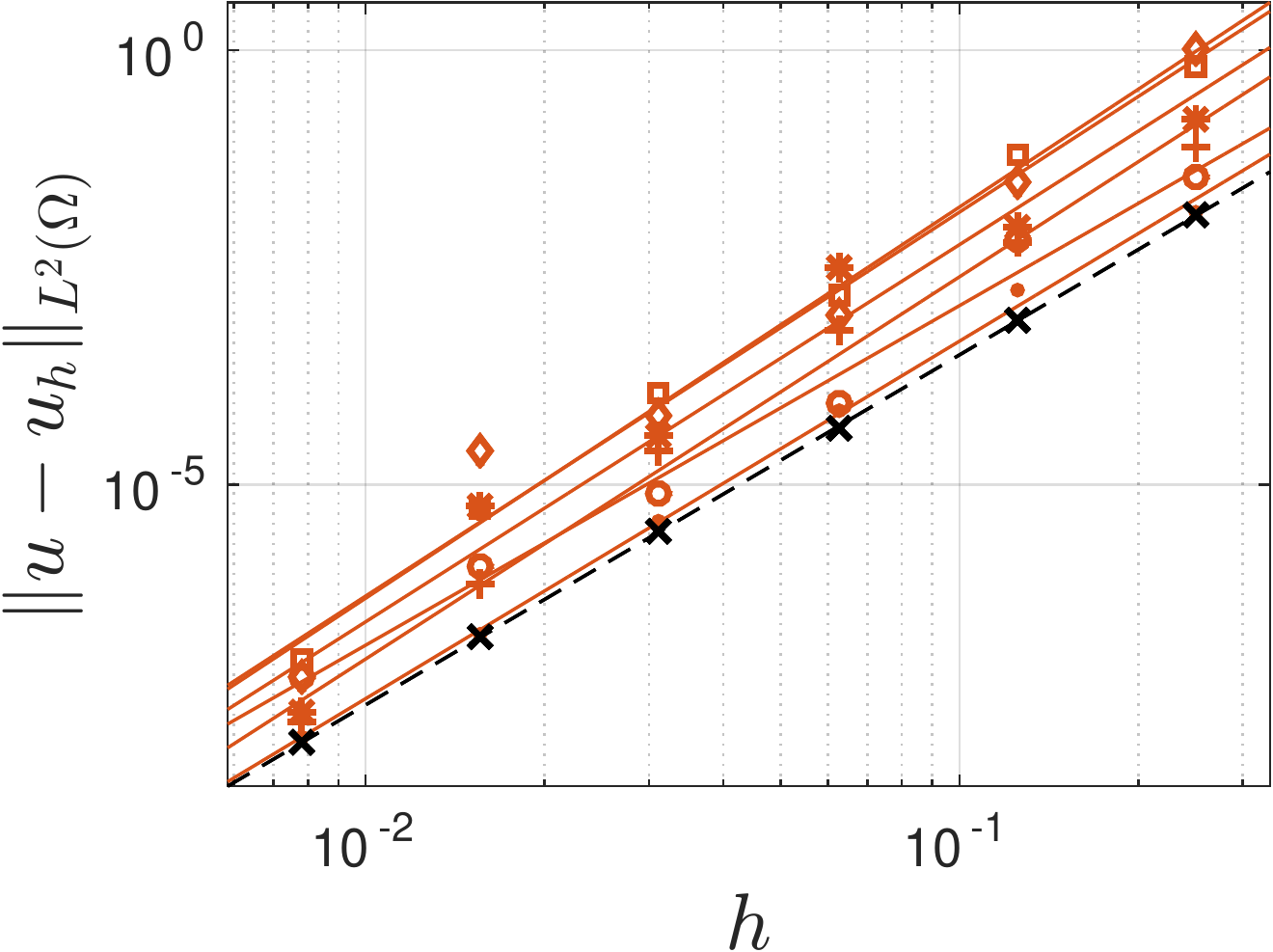}
  \includegraphics[width=0.32\textwidth]{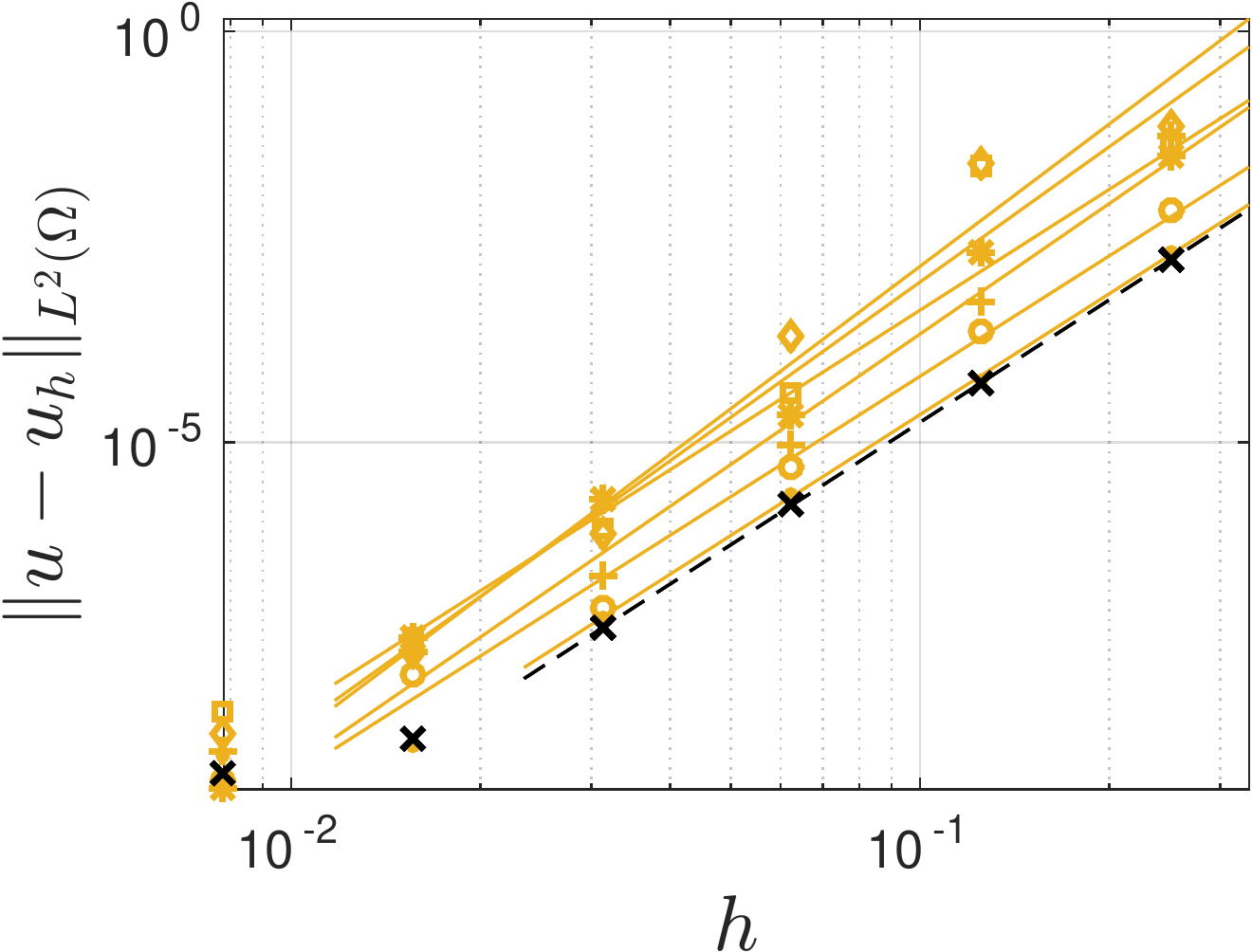}

  \vspace{0.5cm}

  \includegraphics[width=0.32\textwidth]{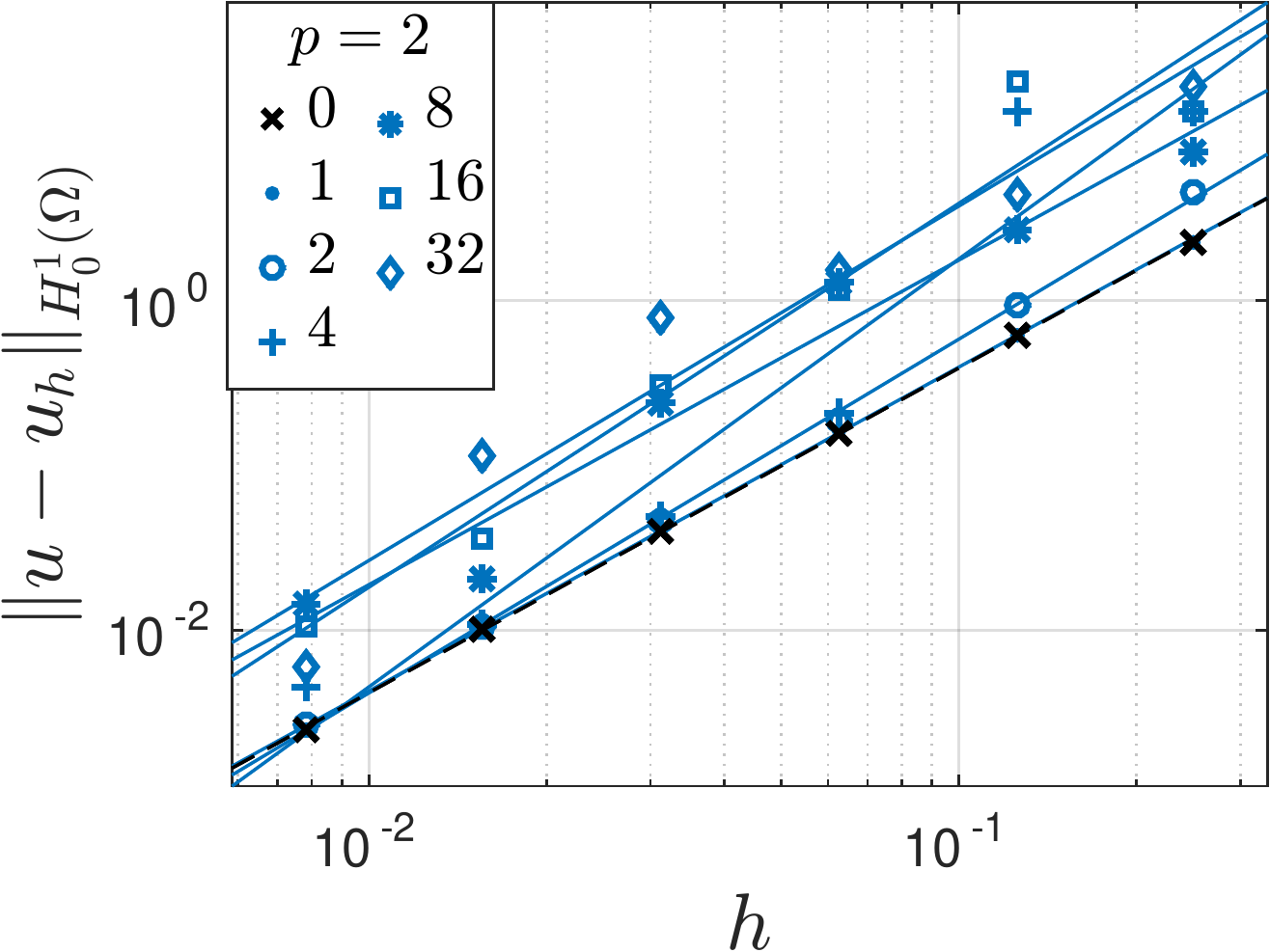}
  \includegraphics[width=0.32\textwidth]{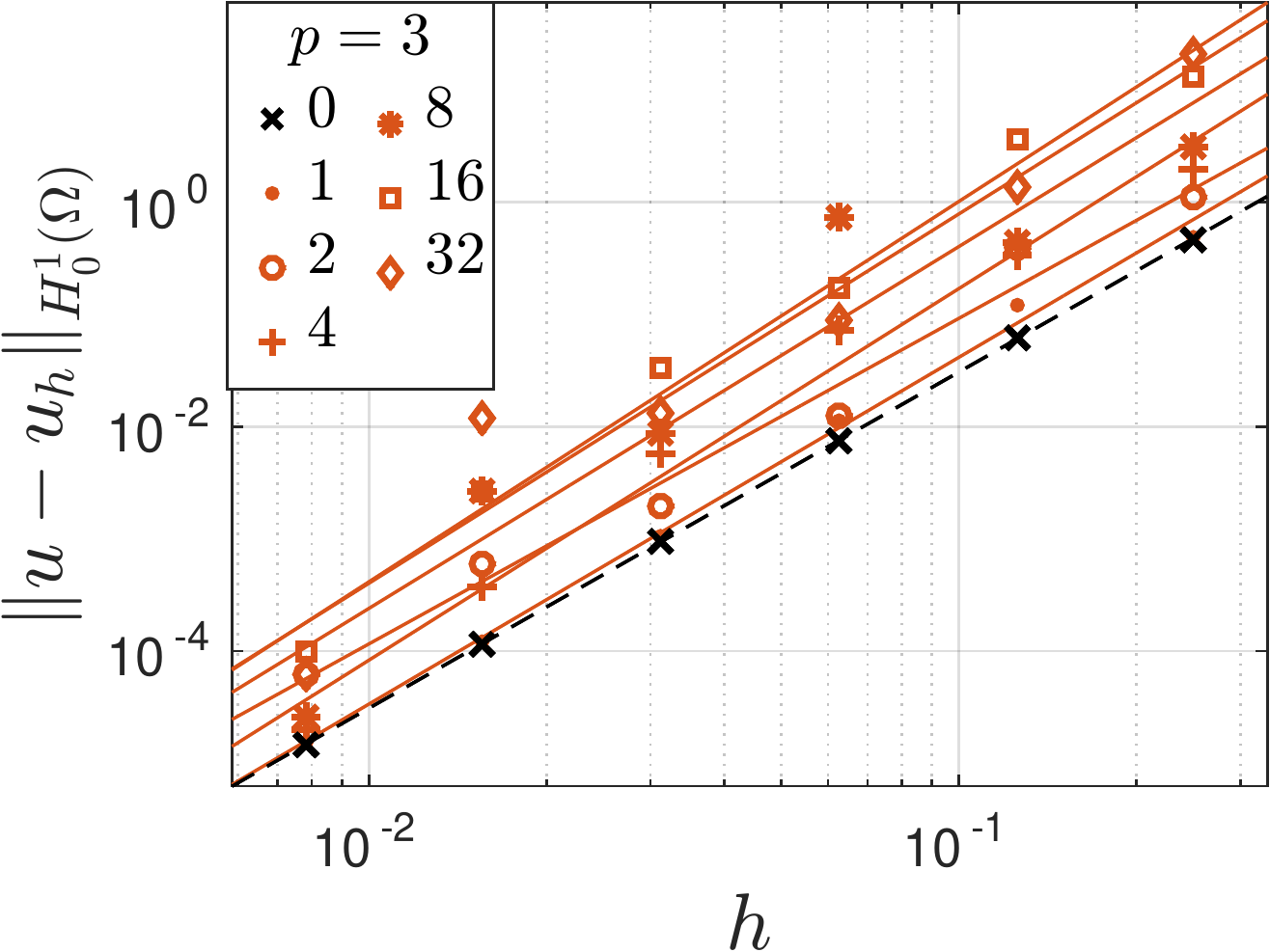}
  \includegraphics[width=0.32\textwidth]{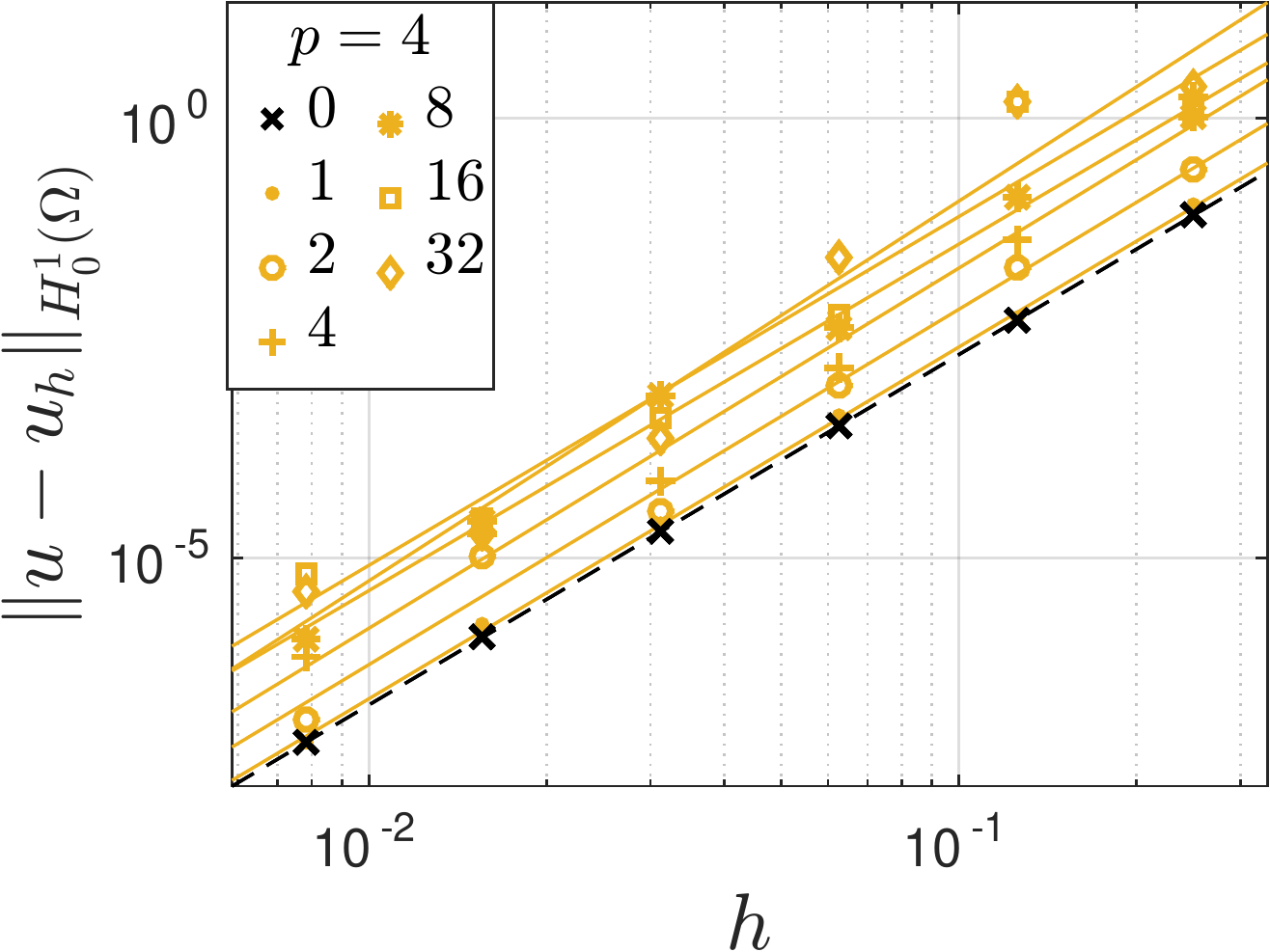}

  \vspace{0.5cm}

  \includegraphics[width=0.32\textwidth]{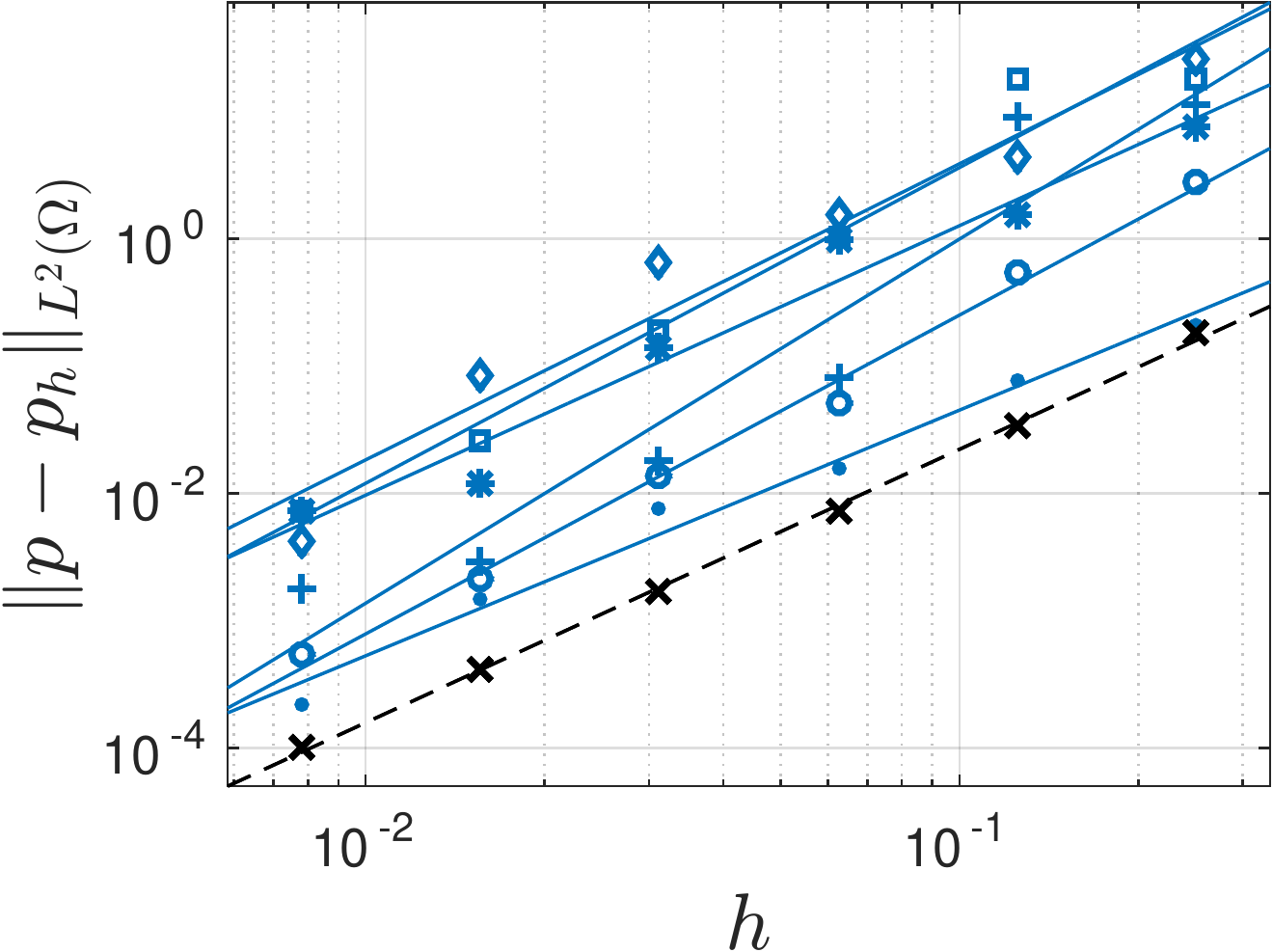}
  \includegraphics[width=0.32\textwidth]{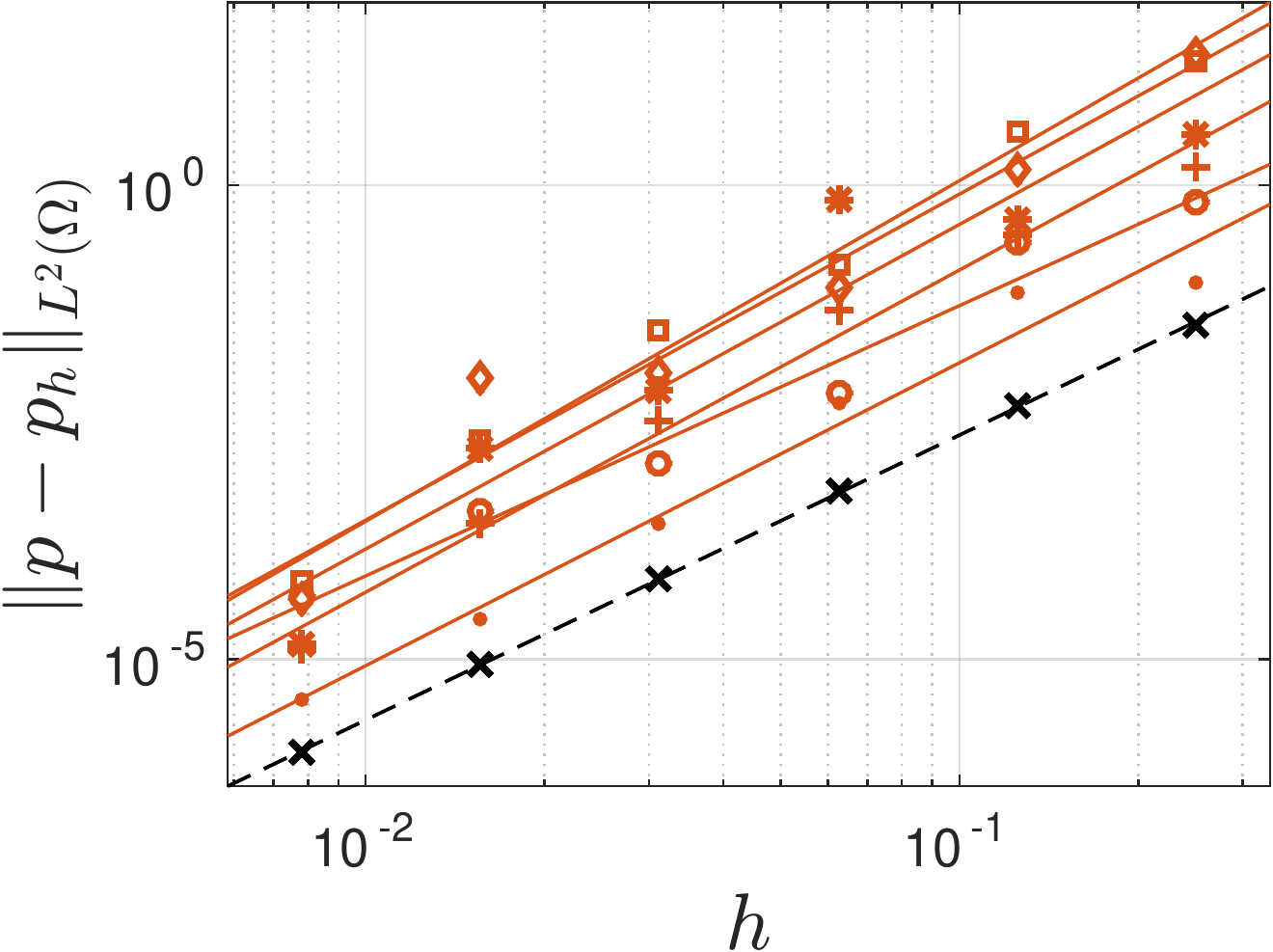}
  \includegraphics[width=0.32\textwidth]{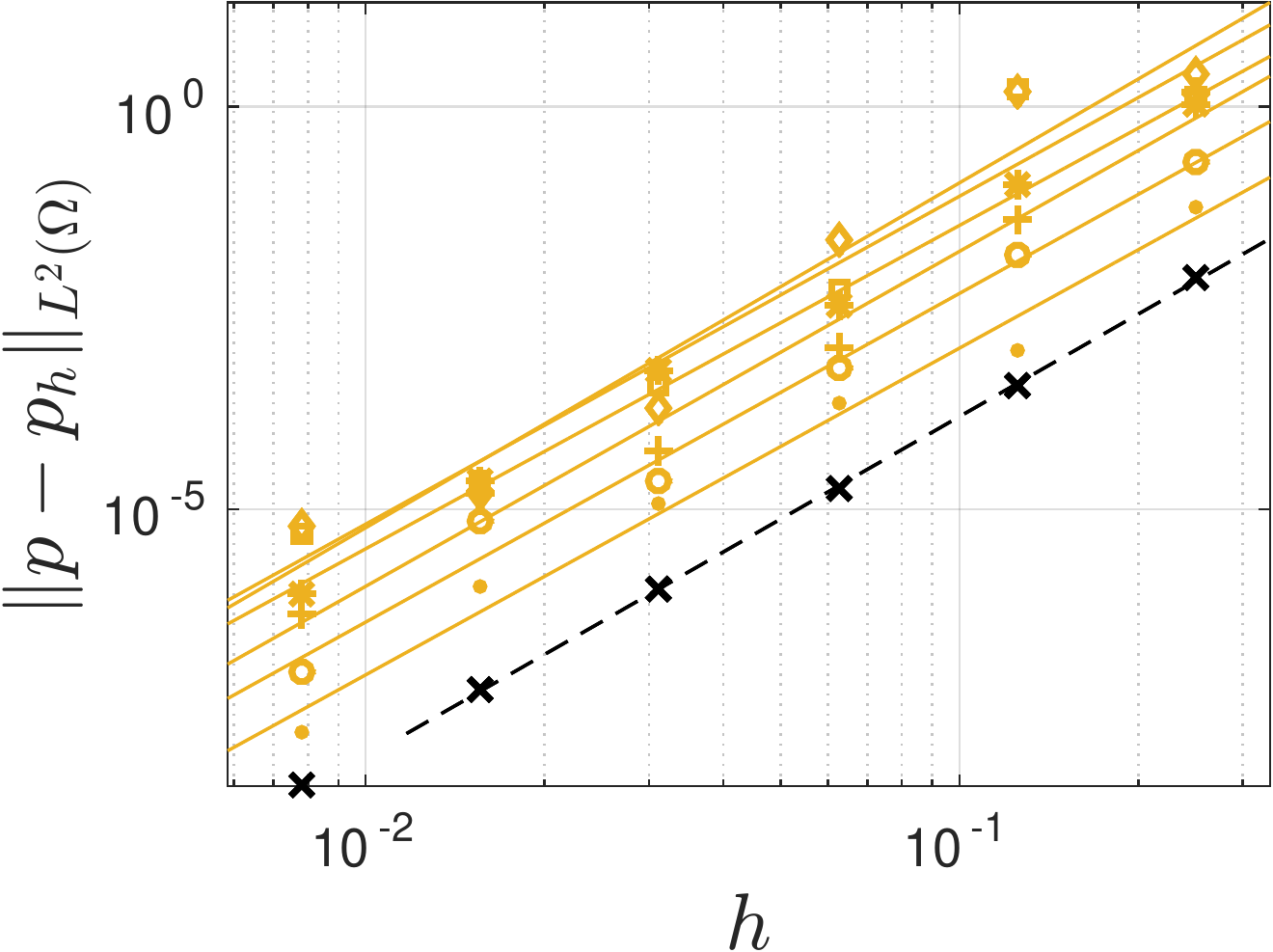}
  \caption{Convergence results for $k=2$, $3$ and $4$ (left to right) using up to $32$ meshes (single mesh results are $N=0$). From top to bottom we have the velocity error in the $L^2(\Omega)$ norm, the velocity error in the $H^1_0(\Omega)$ norm, and the pressure error in the $L^2(\Omega)$ norm. Results less than $10^{-8}$ are not included in the convergence lines due to limits in floating point precision. }
  \label{fig:convplots}
\end{figure}

To investigate the convergence of the multimesh finite element method, we solve the Stokes problem in the unit square with the following exact solution
\begin{align}
  \label{eq:exactu}
  \bfu(x, y) &= 2 \pi \sin(\pi x) \sin(\pi y) \cdot
  ( \cos(\pi y) \sin(\pi x) , -\cos(\pi x) \sin(\pi y) ), \\
  \label{eq:exactp}
  p(x, y) &= \sin(2 \pi x) \sin(2 \pi y),
\end{align}
and corresponding right hand side. We use $P_kP_{k-1}$ Taylor--Hood elements with $k\in \{2,3,4\}$ and we use $N\in \{1,2,4,8,16,32\}$ randomly placed domains as in \cite{mmfem-1} (see Figure~\ref{fig:poisson_meshes}). Due to the random placement of domains, some domains are completely hidden and will not contribute to the solution. For $N=8$, this is the case for one domain, for $N=16$, three domains and for $N=32$, four domains are completely hidden. This is automatically handled by the computational geometry routines. Convergence results are presented in Figure~\ref{fig:convplots} as well as in Table~\ref{table:rates}.

\begin{table}
  \centering
  \caption{Error rates for $\bfe_{L^2} = \|\bfu - \bfu_h\|_{L^2(\Omega)}$, $\bfe_{H^1_0} = \|\bfu - \bfu_h\|_{H^1_0(\Omega)}$ and $e_{L^2} = \|p - p_h\|_{L^2(\Omega)}$.}
  \label{table:rates}
  \begin{tabular}{|c | c c c | c c c | c c c|}
    \toprule
    & & $k=2$  & & & $k=3$  &  & & $k=4$  &   \\
    \midrule
    $N$&
    $\bfe_{L^2}$ & $\bfe_{H^1_0}$ & $e_{L^2}$ &
    $\bfe_{L^2}$ & $\bfe_{H^1_0}$ & $e_{L^2}$ &
    $\bfe_{L^2}$ & $\bfe_{H^1_0}$ & $e_{L^2}$ \\
    \midrule
    0 & 2.9952 & 1.9709 & 2.1466 & 4.0289 & 2.9966 & 3.0028 & 4.9508 & 3.9844 & 4.2286 \\
    1 & 2.9750 & 1.9658 & 1.9291 & 4.1153 & 3.0912 & 3.1932 & 4.8861 & 4.0006 & 4.0587 \\
    2 & 3.2764 & 2.1472 & 2.5036 & 3.9087 & 2.9021 & 2.8489 & 4.8677 & 4.0416 & 4.0832 \\
    4 & 3.6666 & 2.5971 & 2.8597 & 4.3996 & 3.3125 & 3.3957 & 5.2741 & 4.0966 & 4.1609 \\
    8 & 3.0359 & 1.9697 & 2.1163 & 4.3412 & 3.2258 & 3.4213 & 4.8840 & 3.9409 & 4.0169 \\
    16 & 3.4131 & 2.3298 & 2.4794 & 4.5033 & 3.3907 & 3.5910 & 5.4702 & 3.9664 & 4.0729 \\
    32 & 3.2832 & 2.1505 & 2.3255 & 4.4196 & 3.2922 & 3.4362 & 5.7538 & 4.3191 & 4.2848 \\
    \bottomrule
  \end{tabular}
\end{table}

%%%%%%%%%%%%%%%%%%%%%%%%%%%%%%%%%%%%%%%%%%%%%%%%%%%%%%%%%%%%%%%%%%%%%%%%%%%%%%%%
\section{Discussion}

The results presented in Table~\ref{table:rates} and Figure~\ref{fig:convplots} show the expected order of convergence for the velocity in the $L^2(\Omega)$ norm ($k + 1$), for the velocity in the $H^1_0(\Omega)$ norm ($k$), and for the pressure in the $L^2(\Omega)$ norm ($k$).

A detailed inspection of Figure~\ref{fig:convplots} reveals that, as expected, the multimesh discretization yields larger errors than the single mesh discretization (standard Taylor--Hood on one single mesh). The errors introduced by the multimesh discretization are one to two orders of magnitude larger than the single mesh error. However, the convergence rate is optimal and it should be noted that the results presented here are for an extreme scenario where a large number of meshes are simultaneously overlapping; see Figure~\ref{fig:poisson_meshes}. For a normal application, such as the simulation of flow around a collection of objects, each object would be embedded in a boundary-fitted mesh and only a small number of meshes would simultaneously overlap (in addition to each mesh overlapping the fixed background mesh), corresponding to the situation when two or more objects are close.

The presented method and implementation demonstrate the viability of the multimesh method as an attractive alternative to existing methods for discretization of PDEs on domains undergoing large deformations. In particular, the discretization and the implementation are robust to thin intersections and rounding errors, both of which are bound to appear in a simulation involving a large number of meshes, timesteps or configurations.

%%%%%%%%%%%%%%%%%%%%%%%%%%%%%%%%%%%%%%%%%%%%%%%%%%%%%%%%%%%%%%%%%%%%%%%%%%%%%%%%
\begin{acknowledgement}
  August Johansson was supported by the Research Council of Norway through the FRIPRO Program at Simula Research Laboratory, project number 25123. Mats G.\ Larson was supported in part by the Swedish Foundation for Strategic Research Grant No.\ AM13-0029, the Swedish Research Council Grants Nos.\  2013-4708, 2017-03911, and the Swedish Research Programme Essence. Anders Logg was supported by the Swedish Research Council Grant No.\ 2014-6093.
\end{acknowledgement}

%%%%%%%%%%%%%%%%%%%%%%%%%%%%%%%%%%%%%%%%%%%%%%%%%%%%%%%%%%%%%%%%%%%%%%%%%%%%%%%%
\bibliographystyle{spmpsci}
\bibliography{bibliography}
\end{document}